\documentclass{article}
\usepackage{amsfonts}
\usepackage{amssymb}

\begin{document}
\newtheorem{theor}{Theorem}[section] 
\newtheorem{prop}[theor]{Proposition} 
\newtheorem{cor}[theor]{Corollary}
\newtheorem{lemma}[theor]{Lemma}
\newtheorem{sublem}[theor]{Sublemma}
\newtheorem{defin}[theor]{Definition}
\newtheorem{conj}[theor]{Conjecture}

\hfuzz2cm

\gdef\beginProof{\par{\bf Proof: }}
\gdef\endProof{${\bf Q.E.D.}$\par}  
\gdef\ar#1{\widehat{#1}}
\gdef\pr{^{\prime}}
\gdef\prpr{^{\prime\prime}}
\gdef\mtr#1{\overline{#1}}
\gdef\ra{\rightarrow}
\gdef\Bbb{\bf }
\gdef\P1{{\Bbb P}^{1}_{D}}
\gdef\dbd{dd^{c}}
\gdef\a{\alpha}
\gdef\ca{ch_{g}(\alpha)}
\gdef\ttmk#1{\widetilde{\theta}^{k}(#1)}
\gdef\tdm#1{\theta^{k}({#1}^{\vee})^{-1}}
\gdef\td#1{\theta^{k}({#1}^{\vee})}
\gdef\rl{{\Lambda}}
\def\covol{{\rm covol}}
\def\lcovol{{\rm lcovol}}
\gdef\CT{CT}
\gdef\refeq#1{(\ref{#1})}
\gdef\mn{{\mu_{n}}}
\gdef\zn{{\Bbb Z}/(n)}
\gdef\umn{^{\mn}}
\gdef\lmn{_{\mn}}
\gdef\blb{{\big(}}
\gdef\brb{{\big)}}
\gdef\ttmk#1{\widetilde{\lambda}(#1)}
\gdef\chge#1{ch_{g}^{-1}(\lambda_{-1}({#1}^{\vee}))
ch_{g}(\lambda_{-1}({#1}^{\vee}_{g}))}
\gdef\mlcovol{{\rm mlc}}
\gdef\Hom{{\rm Hom}}
\def\Trs{{\rm Tr_s\,}}
\def\Tr{{\rm Tr\,}}
\def\End{{\rm End}}
\def\eq{equivariant }
\def\Td{{\rm Td}}
\def\ch{{\rm ch}}
\def\torus{{\cal T}}
\def\Proj{{\rm Proj}}
\def\bmn{{R_{n}}}
\def\uexp#1{{{\rm e}^{#1}}}
\def\Spec{{\rm Spec}\,}
\def\Qb{\mtr{\Bbb Q}}
\def\Cn{{\bf C}_n}
\def\Zn{{\Bbb Z}/n}
\def\deg{{\rm deg}\,}
\def\mod{{\rm mod}}
\def\ac1{\ar{\rm c}_{1}}
\def\boxtimes{{\otimes_{\rm Ext}}}
\def\Qmm{{{\Bbb Q}(\mu_{m})}}
\def\NIm#1{{\rm Im}(#1)}
\def\NRe#1{{\rm Re}(#1)}
\def\rk{{\rm rk}\,}

\author{ Kai K\"ohler
\footnote{Mathematisches Institut, Wegelerstr. 10, 
D-53115 Bonn, Germany, E-mail: koehler@rhein.iam.uni-bonn.de, 
URL: http://www.math.uni-bonn.de/people/koehler}\\
         Damian Roessler
\footnote{Centre de Math\'ematiques de Jussieu,  
     Universit\'e Paris 7 Denis Diderot,  
     Case Postale 7012,  
     2, place Jussieu,  
     F-75251 Paris Cedex 05, France,   
     E-mail : roessler@math.jussieu.fr,  
     URL: http://www.math.jussieu.fr/$\tilde{\ }$roessler}\\}

\title{A fixed point formula of Lefschetz type in Arakelov geometry II: 
a residue formula / Une formule du point fixe de type Lefschetz 
en g\'eom\'etrie d'Arakelov II: une formule des r\'esidus}
\maketitle
\begin{abstract}
This is the second of a series of papers dealing with an analog 
in Arakelov geometry of the holomorphic Lefschetz fixed point formula. We use 
the 
main result \cite[Th. 4.4]{KR2} of the first paper to prove a 
residue formula "\`a la Bott" for arithmetic characteristic classes living on 
arithmetic varieties acted upon by 
a diagonalisable torus; recent results of Bismut-Goette 
on the equivariant (Ray-Singer) analytic torsion play a key role in 
the proof.\ \ /\ \  
Cet article est le second d'une s\'erie d'articles dont l'objet est 
un analogue en g\'eom\'etrie d'Arakelov 
de la formule du point fixe de Lefschetz holomorphe. Nous utilisons 
le r\'esultat principal \cite[Th. 4.4]{KR2} du premier article pour prouver 
une formule des r\'esidus "\`a la Bott" pour des classes caract\'eristiques 
vivant sur des vari\'et\'es arithm\'etiques munis d'une action de tore; 
de r\'ecents r\'esultats de Bismut-Goette sur la torsion analytique 
\'equivariante (de Ray-Singer) joue un r\^ole cl\'e dans la preuve. 
\end{abstract}
\begin{center}
\end{center}
\begin{center}
\end{center}
\thispagestyle{empty}
\newpage
\setcounter{page}{1}
\tableofcontents
\newpage
\parindent=0pt
\parskip=5pt

\section{Introduction}

This is the second of a series of four papers on 
equivariant Arakelov theory and a fixed point formula therein. 
We give here an application of the main 
result \cite[Th. 4.4]{KR2} of the first paper (announced in \cite{KR1}).
\\
We prove 
 a residue formula "\`a la Bott" (Theorem \ref{residue}) for  the arithmetic
Chern numbers of arithmetic varieties endowed with 
the action of a diagonalisable torus. More precisely, this formula computes 
arithmetic Chern numbers of equivariant Hermitian vector bundles 
(in particular, the height relatively to some equivariant 
projective embedding) as a contribution of arithmetic Chern numbers 
of bundles living of the fixed scheme and an anomaly term, which depends 
on the complex points of the variety only.
Our determination of the anomaly term relies heavily on 
recent results by Bismut-Goette (\cite{BG}, \cite{BG1}). The formula 
in \ref{residue} is formally similar to Bott's residue 
formula \cite[III, Prop. 8.13, p. 598]{AS} for the characteristic numbers of 
vector bundles, up 
to the anomaly term.  
Our method of proof is similar to
Atiyah-Singer's and is described in more detail in the introduction to 
section 2. The effective computability of the anomaly term 
is also discussed there. 
\\
Apart from the residue formula itself, this article has the following 
two side results, which are of independent interest and which we choose to
highlight here, lest they remain unnoticed in the body of the proof of 
Th. \ref{residue}. 
 The first one is a 
corollary of the 
residue formula, which shows that the height relatively to equivariant 
line bundles on torus-equivariant arithmetic varieties depends on 
less data than on general varieties (see corollary \ref{impcor}):\\
{\bf Proposition}. {\it 
Let $Y$ be an arithmetic variety endowed with a torus action. Write 
$Y_\torus$ for the fixed point scheme of $Y$. Suppose 
that $\mtr{L}, \mtr{L}\pr$ are torus-equivariant hermitian line bundles.
 If there is an equivariant isometry $\mtr{L}_{Y_\torus}\simeq \mtr{L}\pr_{Y_\torus}$ over 
 $Y_\torus$ and an equivariant (holomorphic) isometry 
 $\mtr{L}_{\bf C}\simeq\mtr{L}\pr_{\bf C}$ over $Y_{\bf C}$ then 
 the height of $Y$ relatively to $\mtr{L}$ is equal to the height 
 of $Y$ relatively to $\mtr{L}\pr$.}\\
The second one is a conjecture which naturally arises in the course 
of the proof of the residue formula (see lemma \ref{rational}):\\
{\bf Conjecture.} {\it Let $M$ be a $S^{1}$-equivariant projective complex 
manifold, equipped with an $S^{1}$-invariant K\"ahler metric. 
Let $\mtr{E}$ be a $S^{1}$-equivariant complex vector bundle 
on $M$, equipped with a $S^{1}$-invariant hermitian metric. Let 
$T_{g_{t}}(\cdot)$ (resp. $R_{g_t}(\cdot)$ of $E$) be the equivariant analytic torsion 
of $\mtr{E}$ (resp. the equivariant $R$-genus), 
 with respect to the automorphism $e^{2i\pi t}$. 
There is a rational function $Q$ with complex coefficients and 
a pointed neighborhood $\stackrel{\circ}{U}$ of 
$0$ in $\bf R$ such that
$$
{1\over 2}T_{g_t}(M,\mtr{E})-
{1\over 2}\int_{M_{g_t}}\Td_{g_t}(TM)
\ch_{g_t}(E)R_{g_t}(TM)=Q(e^{2\pi it})
$$
if $t\in \stackrel{\circ}{U}$}\\
(here $M_{g_t}$ is the fixed point set of the automorphism $e^{2i\pi t}$, 
$\ch_{g_{t}}$ is the equivariant Chern character and 
$\Td_{g_t}$ is the equivariant Todd genus - see section 4 of 
\cite{KR2} for more details).\\
The lemma \ref{rational} shows that this conjecture is verified, when 
the geometric objects appearing in it have certain models over the integers but 
it seems unlikely that the truth of the conjecture should be 
dependent on the existence of such models.\\
The appendix is logically independent of the rest of the article. We 
formulate a conjectural generalisation of the main result of 
\cite{KR2}.
\\ 
The notations and conventions of the section 4 of \cite{KR2}  
(describing the main result) and 6.2 (containing 
a translation of the fixed point formula into arithmetic Chow theory) 
will be used without comment.  This
article is a part of the habilitation thesis of the first author.
\\
{\bf Acknowledgments.} It is a pleasure to thank Jean-Michel Bismut,
 Sebastian Goette, Christophe 
Soul\'e and Harry Tamvakis for stimulating discussions and 
interesting hints. We thank the SFB 256, "Nonlinear Partial Differential
Equations", at the University of Bonn for its support. The second author 
is grateful to the IHES (Bures-sur-Yvette) and its very able staff 
for its support. 

\section{An "arithmetic" residue formula}

In this subsection, we consider 
arithmetic varieties endowed with an action of 
a diagonalisable torus. We shall use the fixed point formula 
\cite[Th. 4.4]{KR2} 
 to obtain a formula computing 
arithmetic characteristic numbers (like the height 
relatively to a Hermitian line bundle) in terms 
of arithmetic characteristic numbers of the fixed point scheme 
(a "residual" term) and an anomaly 
term derived from the equivariant and non-equivariant 
analytic torsion. One can express this term using
characteristic currents only, without involving the analytic torsion (see 
subsection
2).   See equation \refeq{ResFV} for a first version of the residue formula 
(where
the anomaly  term is expressed via the analytic torsion) and \ref{residue}  for 
the
final formula (where  the anomaly term  is expressed using a characteristic 
current). 
One can use the residue formula to compute the height of 
some flag varieties; 
there the anomaly term can be computed using the explicit values for the torsion 
given
in
\cite{Koehler2}. We shall nevertheless 
not carry out the details of this application, as the next 
paper \cite{KK} gives a general formula for the height of 
flag varieties.  \\
The strategy of proof we follow here is parallel to Atiyah-Singer's 
in \cite[Section 8]{AS}. Notice however that our proof, which 
involves the $\gamma$-operations, works in the algebraic case as well. 
 The fundamental step of the proof is a passage to the limit 
on both sides of the arithmetic fixed point formula, where 
the limit is taken on finite group schemes of increasing order 
inside a given torus. Both sides of the fixed point formula 
can be seen as rational functions of a circle element near $1$
and one can thus identify their constant coefficients. The constant 
coefficient of the arithmetic Lefschetz trace is the arithmetic 
Euler characteristic, which can in turn be related with 
arithmetic characteristic numbers via the (arithmetic) Riemann-Roch 
formula.\\
Furthermore, following a remark of J.-M. Bismut, we would like to point out that 
a direct proof of 
the formula 
\ref{residue} seems tractable. One could proceed as in the proof 
of the fixed point formula \cite[Th. 4.4]{KR2} (by deformation to 
the normal cone) and replace at each step the anomaly formulae for 
the equivariant analytic torsion by the anomaly formulae for 
the integral appearing in \ref{residue}, the latter formulae having much 
easier proofs (as they do not involve the spectrum of Laplace operators). 
One would thus avoid mentioning the analytic torsion altogether. If 
\cite[Th. 4.4]{KR2} and the work of Bismut-Goette was not available, this 
would probably be the most natural way to approach the residue formula.
\\
Let $\torus:=\Spec\ {\Bbb Z}[X,X^{-1}]$ be the one-dimensional 
torus over $\Bbb Z$. 
Let $f:Y\ra \Spec{\Bbb Z}$ be a regular scheme, flat over $\Bbb Z$, 
endowed with a $\torus$-projective action and 
such that the fixed scheme $Y_{\torus}$ is flat over $\Bbb Z$ 
(this requirement is only necessary because we choose 
to work with arithmetic Chow theory).  Let $d+1$ be the absolute 
dimension of $Y$. 
This action induces a holomorphic group action of the multiplicative 
group 
${\Bbb C}^{*}$ on the manifold $Y({\Bbb C})=:M$ and thus an action of 
the circle $S^{1}\subseteq {\Bbb C}^{*}$.  
We equip $Y({\Bbb C})$ once and for all with an 
$S^{1}$-invariant K\"ahler 
metric ${\omega}^{TY({\Bbb C})}=\omega^{TM}$ (such a metric can be obtained  
explicitly via an embedding into some projective space).   
Now let $m>0$ be a strictly positive integer coprime to $n$. 
Consider the homomorphism $s_{m,n}:{\Bbb Z}\ra\Zn$, given by the formula 
$a\mapsto m.(a\,\mod\,n)$. This homomorphism induces an immersion 
$i_{m,n}:\mn\ra\torus$ of group schemes. Let now $E$ be a $\torus$-equivariant 
bundle on $Y$. 
Recall that the equivariant structure of $E$ induces a $\Bbb Z$-grading 
on the restriction $E|_{Y_{\torus}}$ of $E$ to the fixed point scheme 
of the action of $\torus$ on $Y$; the $k$-th term ($k\in{\Bbb Z}$) of this 
grading 
is then denoted by $E_{k}$. 
\begin{lemma}
Write $E^{m,n}$ for $E$ viewed as a $\mn$-equivariant 
bundle via $i_{m,n}$. 
There exists an $\epsilon>0$ such that 
for all $k\in{\Bbb Z}$ 
the natural injection $E_{k}\ra E_{s_{m,n}(k)}^{n,m}$ is 
an isomorphism if $1/n<\epsilon$. 
\label{ContLem}
\end{lemma}
\beginProof
This natural injection is an isomorphism iff the equality 
$s_{m,n}(k)=s_{m,n}(k\pr)$ 
($k,k\pr\in{\Bbb Z}$) implies that $k=k\pr$. Now notice 
that the kernel of $s_{m,n}$ is generated by $n$. 
Thus the implication is realized if we choose $\epsilon$ 
such that $1/\epsilon>2.\max\{|k|\,|k\in{\Bbb Z}, E_{k}\not= 0\}$ and 
we are done. 
\endProof
\begin{cor}
Let $P$ be a projective space over $\Bbb Z$ endowed 
with a global action of the torus $\torus$. Write $P^{m,n}$ for $P$ 
viewed as a $\mn$-equivariant 
scheme via $i_{m,n}$. Then there 
exists $\epsilon>0$, such that if $1/n<\epsilon$, then 
the closed immersion 
$P_{\cal T}\ra P^{m,n}_{\cal T}$ is an isomorphism.
\end{cor}
\beginProof Let $M$ be a free $\bf Z$-module endowed 
with a $\torus$-action, such that there is an equivariant 
isomorphism $P\simeq{\rm\bf P}(M)$. Let us write $M^{m,n}$ for $M$ viewed as a 
$\mn$-comodule via $i_{m,n}$. 
By the  
description of the fixed scheme given in \cite[Prop. 2.9]{KR2}, we have 
$P_{\torus}=\coprod_{k\in {\bf Z}}{\rm\bf P}(M_{k})$ and 
$P_{\torus}^{m,n}=\coprod_{k\in \Zn}{\rm\bf P}(M_{k}^{m,n})$. 
Furthermore, by construction the immersion 
${\rm\bf P}(M_{k})\ra P$ factors through 
the immersion ${\rm\bf P}(M_{k})\ra{\rm\bf P}(M_{s_{m,n}(k)}^{m,n})$ 
induced by the injection $M_{k}\ra M_{s_{m,n}(k)}^{n,m}$. 
By the last lemma, there exists an $\epsilon>0$ such that 
for all $k\in{\Bbb Z}$ 
the natural injection $M_{k}\ra M_{s_{m,n}(k)}^{n,m}$ is 
an isomorphism if $1/n<\epsilon$. From this, we can conclude. 
\endProof 
Let again $E$ be a $\torus$-equivariant bundle on 
$Y$, such that the cohomology of $E$ vanishes 
in positive degrees. We equip $E_{\Bbb C}$ with an
$S^{1}$-invariant  metric (such a metric can be obtained from an 
arbitrary metric by integration). 
Consider $E$ and $Y$ as $\mu_{n}$-equivariant via $i_{m,n}$. 
We shall apply \cite[Th. 6.14]{KR2} to $E$. For this application, 
we fix 
the primitive root of unity $\uexp{2i\pi m/n}$ of $\mn({\Bbb C})$. 
If $\alpha\in {\bf C}^*$, we shall write $g(\alpha)$ for the 
corresponding automorphism of $Y({\Bbb C})$ and we let
$g_{m,n}:=g(\uexp{2i\pi m/n})$.  Set $M:=Y({\Bbb C})$.
By \cite[Th. 6.14]{KR2}, we get 
\begin{eqnarray}
\ar{\deg}_{\mn}(R^{0}f_{*}(\mtr{E}))&=&{1\over
2}T_{{{g_{m,n}}}}(M,\mtr{E}) -
{1\over 2}\int_{M_{g_{m,n}}}\Td_{{{g_{m,n}}}}(TM)
\ch_{{g_{m,n}}}(E)R_{{{g_{m,n}}}}(TM)
\nonumber
\\&&+
\ar{{\rm
deg}}\blb f_*^{\mn}\blb\sum_{i=0}^{\rk(N^{\vee}_{Y/Y_{\mn}})}
(-1)^{i}\ar{\ch}_{\mn}(\Lambda^{i}(\mtr{N}^{\vee}
_{Y/Y_{\mn}}))\brb^{-1}
\nonumber\\&&.
\ar{\Td}(\mtr{Tf^{\mn}}).
\ar{\ch}_{\mn}(\mtr{E})\brb
\label{reform}
\end{eqnarray}
Furthermore, 
using the last lemma and its corollary, we see that there is 
an $\epsilon>0$ and a rational function ${\cal Q}(\cdot)$ with 
coefficients in $\ar{\rm CH}(Y_{\torus})_{\bf C}$, such that for all $n,m$ 
coprime with  
$1/n<\epsilon$, the term 
\begin{equation}
 \blb\sum_{i=0}^{\rk(N^{\vee}_{Y/Y_{\mu_{n}}})}
(-1)^{i}\ar{\ch}_{\mu_{n}}(\Lambda^{i}(\mtr{N}^{\vee}
_{Y/Y_{\mu_{n}}}))\brb^{-1}
.\ar{\Td}(\mtr{Tf^{\mn}}).
\ar{\ch}_{\mn}(\mtr{E})
\label{suplabel1}
\end{equation}
equals ${\cal Q}(\uexp{2i\pi m/n})$. Similarly , 
there is 
an $\epsilon>0$ and a rational function ${Q}(\cdot)$ with 
complex coefficients, such that for all $n,m$ coprime with  
$1/n<\epsilon$, the term 
\begin{equation}
\ar{{\rm deg}}\blb
f_*^{\mu_{n}} \blb\sum_{i=0}^{\rk(N^{\vee}_{Y/Y_{\mu_{n}}})}
(-1)^{i}\ar{\ch}_{\mu_{n}}(\Lambda^{i}(\mtr{N}^{\vee}
_{Y/Y_{\mu_{n}}}))\brb^{-1}
.\ar{\Td}(\mtr{Tf^{\mn}}).
\ar{\ch}_{\mn}(\mtr{E})\brb
\label{suplabel}
\end{equation}
equals $Q(\uexp{2i\pi m/n})$. 
Since the elements of the type $\uexp{2i\pi t}$, where $t\in{\Bbb Q}$, 
form a dense subset of $S^{1}$, we see that the function $Q(z)$ is uniquely 
determined. 
Let us call 
$A_{\torus}(\mtr{E})$ the constant term in the Laurent 
development of $Q(z)$ around $1$. 
By construction, there is a polynomial $P(z)$ with 
complex coefficients, such that 
$
\ar{\deg}_{\mn}(R^{0}f_{*}(\mtr{E}))
$ 
equals $P(\uexp{2i\pi m/n})$. By density again, this polynomial 
is uniquely determined. The constant term of its 
Laurent development around $1$ (i.e. its value at $1$) 
is the quantity $\ar{\deg}(R^0f_{*}\mtr{E})$. 
Using \refeq{reform}, we thus see that there is 
a uniquely determined rational function $Q\pr(z)$ with complex 
coefficients and an $\epsilon>0$, such that the quantity
$$
{1\over 2}T_{{{g_{m,n}}}}(M,\mtr{E})-
{1\over 2}\int_{M_{g_{m,n}}}\Td_{{g_{m,n}}}(TM)
\ch_{{g_{m,n}}}(E)R_{{g_{m,n}}}(TM)
$$
equals $Q\pr(\uexp{2i\pi m/n})$ if $1/n<\epsilon$.

Now notice the following. 
Let $I\subset{\bf R}$ be an interval such that  the fixed point set
$M_{g_t}$ does not vary for $t\in I$. Let $g_t:=g(e^{2\pi it})$. Then 
$R_{g_t}$  varies continuously on $I$ (e.g. using
\cite[Remark p. 108]{Koehler2}).

\begin{lemma}\label{rational}
There
is a pointed neighborhood $\stackrel{\circ}{U}$ of 
$0$ in $\bf R$ such that
$$
{1\over 2}T_{g_t}(M,\mtr{E})-
{1\over 2}\int_{M_{g_t}}\Td_{g_t}(TM)
\ch_{g_t}(E)R_{g_t}(TM)=Q\pr(e^{2\pi it})
$$
if $t\in \stackrel{\circ}{U}$. Furthermore this equality holds for all up to 
finitely many values of $e^{2\pi i t}\in S^1$.
\end{lemma}
\beginProof
It remains to prove that the analytic torsion $T_{g_t}(M,\mtr{E})$ is
continuous in $t$ on $\stackrel{\circ}{U}$ (see also \cite{BG}). Let
$I:=]0,\epsilon[$ be an interval on which the fixed point set
$M_{g_t}$ does not vary. Let $M_{g_t}=\bigcup_\mu
M_\mu$ be the decomposition of the fixed point set into connected components of
dimension $\dim M_\mu=:d_\mu$.

Let
$P^\bot$ denote the projection of $\Gamma^\infty(\Lambda^qT^{*0,1}M\otimes E)$ 
on the
orthogonal complement of the kernel of the Kodaira-Laplace operator
$\square_q$ for $0\leq q\leq d$. 
 As shown in Donnelly \cite[Th. 5.1]{Do}, Donnelly and
Patodi \cite[Th. 3.1]{DoP} (see also \cite[Th. 6.11]{BGV}) the trace of the equivariant heat
kernel of
$\square$ for
$u\to 0$ has an asymptotic expansion providing the formula
$$
\sum_q (-1)^{q+1} q\Tr g_t^* e^{-u \square_q}P^\bot \sim
\sum_{\mu}\sum_{k=-d_\mu}^\infty u^k \int_{M_\mu} b^\mu_k(t,x)\,d{\rm vol}_x
$$
where the $b_k(t,x)$ are rational functions in $t$ which are non-singular on 
$I$.
Thus the analytic torsion is given by
\begin{eqnarray*}\lefteqn{
T_{g_t}(M,\mtr{E})=\int_1^\infty \sum_q (-1)^{q+1} q\Tr g_t^* e^{-u 
\square_q}P^\bot
\frac{du}u
}\\&&
+\int_0^1 \left(\sum_q (-1)^{q+1} q\Tr g_t^* e^{-u \square_q}P^\bot-
\sum_{\mu}\sum_{k=-d_\mu}^0
u^k \int_{M_\mu} b^\mu_k(t,x)\,d{\rm vol}_x
\right) \frac{du}u
\\&&
+\sum_{\mu}\sum_{k=-d_\mu}^{-1}
\frac{1}k \int_{M_\mu} b^\mu_k(t,x)\,d{\rm vol}_x
-\Gamma'(1)\sum_{\mu} \int_{M_\mu} b^\mu_0(t,x)\,d{\rm vol}_x\,\,.
\end{eqnarray*}
The integrand of the first term is uniformly bounded (in $t$) by the 
non-equivariant
heat kernel. Hence we see in particular that $T_{g_t}(M,\mtr{E})$ is continuous 
in
$t\in I$. As the equation in the lemma holds on a dense subset of $I$, it holds 
in
$I$ and by symmetry for a pointed neighborhood of $0$.
\endProof
Recall that $d+1$ is the (absolute) dimension of $Y$. 
Consider the vector field $K\in\Gamma(TM)$ such
that $e^{t K}=g_t$ on $M$. In \cite{Koehler} the function $R^{{\rm rot}}$
on
${\bf R}\setminus 2\pi {\bf Z}$ has been defined as
$$
R^{{\rm rot}}(\phi):=\lim_{s\to 0^+}\frac\partial{\partial s} \sum_{k=1}^\infty
\frac{\sin k\phi}{k^s}
$$
(according to Abel's
Lemma the series in this definition converges for Re $s>0$).
\begin{cor} \label{rigid}
Let $D_\mu\subset{\bf Z}$ denote the set of all non-zero eigenvalues
of the action of $K/2\pi$ on $TM_{|M_\mu}$ at the fixed point component $M_\mu$. 
There are rational functions $Q'', Q_{\mu,k,l}$ for $k\in
D_\mu$, $0\leq l\leq d_\mu$ such that for all but finiteley many values of $e^{2\pi i t}$
$$
T_{g_t}(M,\mtr{E})=Q''(e^{2\pi i t})+\sum_\mu \sum_{k\in D_\mu} 
\sum_{l=0}^{d_\mu}
Q_{\mu,k,l}(e^{2\pi i t})\cdot (\frac\partial{\partial t})^l R^{{\rm rot}}(2\pi 
i k
t)\,\,.
$$
The functions $Q_{\mu,k,l}$ depend only on the holomorphic structure of $E$ and 
the
complex structure on $M$.
\end{cor}
\beginProof
For $\zeta=e^{i\phi}\in S^1$, $\zeta\neq 1$, let $L(\zeta,s)$
denote the zeta function defined in
\cite[section 3.3]{KR2} with $L(\zeta,s)=\sum_{k=1}^\infty k^{-s}\zeta^k$ for Re
$s>0$. In
\cite[equation (77)]{Koehler2} it is shown that $L(\zeta,-l)$ is a rational 
function
in $\zeta$ for $l\in{\bf N}_0$. Also by \cite[equation (80)]{Koehler2},
$$
\frac{\partial}{\partial s}_{|s=-l}(L(e^{i\phi},s)-(-1)^l 
L({e^{-i\phi}},s))
=\left(\frac{-i\partial}{\partial\phi}\right)^l 2iR^{\rm rot}(\phi)\,\,.
$$
The corollary follows by the definition of the Bismut $R_g$-class (see 
\cite[Def. 3.6]{KR2}) and lemma
\ref{rational}.
\endProof
{\bf Remark.} One might reasonably conjecture that the 
Lemma \ref{rational} is valid on any compact K\"ahler manifold endowed 
with a holomorphic action of $S^1$.\\

 Let us call $L_{\torus}(\mtr{E})$ the constant
term  in the Laurent development of $Q\pr(z)$ around $z=1$. By lemma 
\ref{rational}
we obtain
$$
\ar{\deg}(R^{0}f_{*}(\mtr{E}))=L_{\torus}(\mtr{E})+A_{\torus}(\mtr{E}).
$$
Since for any $\torus$-equivariant bundle, one can find a
resolution by acyclic (i.e. whose cohomology vanishes in positive 
degrees) $\torus$-equivariant bundles, one can drop the acyclicity statement 
in the last equation. More explicitly, one obtains
\begin{equation}
\sum_{q\geq 0}(-1)^q\blb\ar{\deg}((R^{q}f_{*}(\mtr{E}))_{\rm free})+
\log(\#(R^{q}f_{*}(\mtr{E}))_{\rm 
Tors})\brb=L_{\torus}(\mtr{E})+A_{\torus}(\mtr{E}).
\label{BC1}
\end{equation}
Notice 
that $Q\pr(z)$ and thus $L_{\torus}(\mtr{E})$ depends on the K\"ahler
form $\omega^{TM}$ and 
$\mtr{E_{\Bbb C}}$ only and can thus be computed without reference 
to the finite part of $Y$.\\
In the next subsection, we shall apply the last equation to a specific 
virtual vector bundle, which has the property that its Chern character 
has only a top degree term and compute $A_{\torus}(\cdot)$ in this case. 
We then obtain a first version of the residue formula, which arises 
from the fact that the left hand-side of the last equation is 
also computed by the (non-equivariant) arithmetic Riemann-Roch. 
The following subsection then shows how $L_{\torus}(\cdot)$ can 
be computed using the results of Bismut-Goette \cite{BG}; combining 
the results of that subsection with the first version of the residue 
formula gives our final version \ref{residue}. 

\subsection{Determination of the residual term}
 
Let $\mtr{F}$ be an $\cal T$-equivariant Hermitian vector bundle on $Y$. 
\begin{defin}
The polynomial equivariant arithmetic total Chern class 
$\ar{c}_{t}(\mtr{F})\in\ar{\rm CH}_{\Bbb C}(Y_{\torus})[t]$ is 
defined by the formula
$$
\ar{c}_{t}(\mtr{F}):=\prod_{n\in{\bf Z}}\sum_{p=
0}^{\rk(F_n)}\sum_{j=0}^{p}{\rk(F_{n})-j\choose p-j}\ar{c}_j(\mtr{F}_{n})(2\pi 
in t)^{p-j}
$$
where $i$ is the imaginary constant.
\end{defin}
We can accordingly define the $k$-th polynomial (equivariant, 
arithmetic) Chern class $\ar{c}_{k,t}(\mtr{F})$ 
of $\mtr{F}$ as the part of 
$\ar{c}_{t}(\mtr{F})$ lying in 
$(\ar{\rm CH}(Y_{\torus})_{\Bbb C}[t])^{(k)}$, where 
$(\ar{\rm CH}(Y_{\torus})_{\Bbb C}[t])^{(k)}$ are the homogeneous 
polynomials of weighted degree $k$ (with respect to the grading of $\ar{\rm 
CH}(Y_\torus)_{\Bbb C}$).  
Define now $\Lambda_{t}(\mtr{F})$ as the formal 
power series $\sum_{i\geq 0}\Lambda^{i}(\mtr{F}).t^{i}$. Let 
$\gamma^{q}(\mtr{F})$ be the $q$-th coefficient in the 
formal power series $\Lambda_{t/(1-t)}(\mtr{F})$; this is 
a $\bf Z$-linear combination of equivariant Hermitian bundles. 
We denote by $\ar{ch}_{t}(\mtr{F})$ the polynomial equivariant  Chern
character and by 
$\ar{ch}_{t}^{q}(\mtr{F})$ the component of 
$\ar{ch}_{t}(\mtr{F})$ lying in $(\ar{\rm CH}(Y_{\torus})_{\Bbb C}[t])^{(q)}$
\begin{lemma}\label{chgamma}
The element $\ar{\ch}^{p}_{t}(\gamma^q(\mtr{E}-\rk(E)))$ is equal 
to $\ar{c}_{q,t}(\mtr{E})$ if $p=q$ and vanishes if $p<q$.
\label{cruxlem}
\end{lemma}
\beginProof
It is proved 
in \cite[II, Th. 7.3.4]{GS3} that $\ar{\ch}$, as a map from the arithmetic 
Grothendieck group $\ar{K}_{0}(Y)$ to the arithmetic Chow theory 
$\ar{\rm CH}(Y)$ is a map of $\lambda$-rings, where 
the second ring is endowed with the $\lambda$-ring structure 
arising from its grading. Thus 
$\ar{\ch}^{p}_t(\gamma^q(\mtr{F}-\rk(F)))$ is a polynomial 
in the Chern classes $\ac1(\mtr{E}),\ar{c}_2(\mtr{E}),\dots$ and the variable 
$t$.
By construction, its coefficients only depend on the equivariant structure 
of $E$ restricted to $Y_{\torus}$. We can thus suppose for the time of this 
proof that 
the action of $\torus$ on $Y$ is trivial. 
To identify these coefficients, we consider the analogous expression 
$\ch^{p}_{t}(\gamma^q({F}-\rk(F)))$ with values in the polynomial ring 
${\rm CH}(Y)[t]$, where ${\rm CH}(Y)$ is the algebraic Chow ring. 
By the same token this is a polynomial in the classical Chern classes 
$c_1(E),c_2(E),\dots$ and the variable $t$. 
As the forgetful map $\ar{\rm CH}(Y)\ra {\rm CH}(Y)$ is 
a map of $\lambda$-rings, the coefficients of these polynomials 
are the same. Thus we can apply the algebraic splitting 
principle and suppose that $F=\oplus_{i=1}^j L_{j}$, where 
the $L_{i}$ are equivariant line bundles. We compute 
$\ch^p_{t}(\gamma^q(F-\rk(F)))=\ch^{p}_{t}(\sigma_{q}(L_{1}-1,\dots ,L_{j}-1))=
(\sigma_{q}(\ch_{t}(L_1)-1,\dots ,\ch_{t}(L_{j})-1))^{(p)}$. 
As the term of lowest degree in 
$ch_{t}(L_{i}-1)$ is $c_{1,t}(L_{i})$, which is 
of (total!) degree $1$, the term of lowest degree in the expression after the 
last equality is $\sigma_{q}(c_{1,t}(L_{1}),\dots 
,c_{1,t}(L_{j}))$ which is of degree $q$ and is equal to $c_{q,t}(F)$ and so 
we are done.
\endProof

{\bf Remark.} An equivariant holomorphic vector bundle $E$ splits
at every component $M_\mu$ of the fixed point set into a sum of vector bundles
$\bigoplus E_\theta$ such that $K$ acts on $E_\theta$ as $i\theta\in i{\bf R}$. 
The $E_\theta$
are those $E_{n,\bf C}$ which do not vanish on $M_\mu$.
Equip $E$ with an invariant Hermitian metric. Then the polynomial
equivariant total Chern form
${c}_{tK}(\cdot)$ is  given by the formula
\begin{equation}\label{eqclass}
{c}_{tK}(\mtr{E})_{|M_\mu}=
\det(\frac{-\Omega^E}{2\pi i}+it \Theta^E+{\rm Id})
=\prod_{\theta\in{\bf R}}\sum_{q=0}^{\rk
E_\theta} {c}_q(\mtr{E}_\theta)(1+it\theta)^{\rk E_\theta-q}
\end{equation}
where $\Theta^E$ denotes the action of $K$ on $E$ restricted to $M_\mu$.
Let $N$ be the normal bundle to the fixed point set. Set
$$
(c^{\rm top}_{tK}(\mtr N)^{-1})':=\frac{\partial}{\partial
b}_{|b=0}c_{\rk N}(\frac{-\Omega^N}{2\pi i}+it\Theta+b\, {\rm Id})^{-1}
$$
where $\Theta$ is the action of $K$ on $N$.
Furthermore, let $r$
denote the additive characteristic class which is given by
\begin{eqnarray*}
r_K(L)_{|M_{\mu}}&:=&
-\frac1{c_1(L)+i
\phi_\mu}\left(-2\Gamma'(1)+2\log|\phi_\mu|+\log(1+\frac{c_1(L)}{i 
\phi_\mu})\right)\\
&=&-\sum_{j\geq 0}\frac{(-c_1(L))^j}{(i\phi_\mu)^{j+1}}\left(
-2\Gamma'(1)+2\log|\phi_\mu|-\sum_{k=1}^j\frac1{k}
\right)
\end{eqnarray*}
for $L$ a line bundle acted upon by $K$ with an angle $\phi_\mu\in{\bf R}$ at 
$M_\mu$ 
(i.e. the Lie derivative by $K$ acts as multiplication by $\phi_\mu$).\\ 
In the next proposition, if $\mtr{E}$ is a  
Hermitian 
equivariant bundle, we write 
$\ar\Td_{g_t}(\mtr {Tf})\ar\ch_{g_t}(\mtr{E})$ for the formal 
Laurent power series development in $t$ 
of the function ${\cal Q}(\uexp{2\pi it})$, where 
${\cal Q}(\cdot)$ is the function defined in \refeq{suplabel1}.  
Set $\ar c^{\rm top}_t(\mtr E):=\ar c_{\rk E,t} (\mtr E)$ for any equivariant 
Hermitian vector bundle $\mtr E$. Note that this class is invertible in 
$\ar{\rm CH}(Y_{\torus})_{\Bbb C}[t]$ if $E_{Y_{\torus}}$ has no invariant 
subbundle.
\begin{prop}\label{regul}
Let $q_{1},\dots,q_{k}$ be natural numbers such that 
$\sum_j{q_j}=d+1$. Let $\mtr{E}^1,\dots,\mtr{E}^k$ be $\torus$-equivariant 
Hermitian bundles. 
Set $x:=\prod_{j} \gamma^{q_j}(\mtr E^j-\rk E^j)$.
The expression
\begin{equation}
\ar\Td_{g_t}(\mtr {Tf})\ar\ch_{g_t}(x)
\label{pfff}
\end{equation}
has a formal Taylor series expansion in $t$. Its constant term is given by
\begin{equation}\label{atregul}
\ar c^{\rm top}_t(\mtr N)^{-1}\prod_j \ar c_{q_j,t}(\mtr E^j)
\end{equation}
which is independent of $t$. 
Also for $t\to 0$
\begin{eqnarray*}
\lefteqn{
\Td_{g_t}(TM)R_{g_t}(TM)\ch_{g_t}(x)}
\\&=&\log (t^2)\cdot (c^{\rm top}_K(N)^{-1})'\prod_j c_{q_j,K}(E^j)
\\&&
+\frac{r_K(N)}{c^{\rm top}_K(N)}\prod_j c_{q_j,K}(E^j)+O(t \log t)\,\,.
\end{eqnarray*}
\end{prop}
Note that the first statement implies that 
$\ar \deg f^\torus_*(\ar c^{\rm top}_t(\mtr N)^{-1}\prod_j \ar
c_{q_j,t}(\mtr E^j))=A_\torus(x).$
\beginProof To prove that the first statement holds, we 
consider that by construction, both the expression 
\refeq{atregul} and the constant term of \refeq{pfff} (as a formal Laurent power 
series) 
are universal polynomials in the Chern classes of the terms 
of the grading of $Tf$ and the terms of the grading of $x$. By using 
Grassmannians (more precisely, products 
of Grassmannians) as in the proof of \ref{cruxlem}, we can reduce the problem 
of the determination of these coefficients to the algebraic case and then 
suppose that 
all the relevant bundles split. Thus, without loss of generality, we consider a
vector bundle
$\mtr E:=\bigoplus_\nu
\mtr L_\nu$ which splits into a direct sum of line bundles $\mtr L_\nu$ on which
$\torus$  acts with multiplicity $m_\nu$. Assume now
$m_\nu\neq 0$ for all $\nu$. Set $x_\nu:=\hat c_1(\mtr L_\nu)$. Then
$$
\ar\Td_{g_t}(\mtr E)=\prod\left(1-e^{-2\pi i t m_\nu-x_\nu}\right)^{-1}.
$$
Now
\begin{eqnarray*}
\left(1-e^{-2\pi i t m_\nu-x_\nu}\right)^{-1}
&=&\frac1{2\pi i t m_\nu+x_\nu}+O(1)\\
&=&\sum_{j=0}^{d-\rk N+1}\frac{(-x_\nu)^j}{(2\pi i t m_\nu)^{j+1}}+O(1)
\end{eqnarray*}
as $t\to 0$.
By definition,
\begin{equation}\label{cInv}
\ar c^{\rm top}_{t}(\mtr E)^{-1}=\prod_\nu \frac1{2\pi i t m_\nu+x_\nu}
=\prod_\nu \sum_{j=0}^{d-\rk N+1}\frac{(-x_\nu)^j}{(2\pi i t m_\nu)^{j+1}}
\,\,.
\end{equation}
Thus,
$\ar\Td_{g_t}(\mtr{E})$ has a Laurent expansion of the form
$$
\ar\Td_{g_t}(\mtr{E})=\ar c^{\rm top}_{t}(\mtr E)^{-1}+\sum_{j=0}^{d-\rk N+1}
t^{-\rk N-j}\hat q_j(t)
$$
with classes $\hat q_j$ of degree $j$ which have a Taylor expansion in $t$ and 
$\hat
q_j(0)=0$. 
As $\ar\Td_{g_t}(\mtr{Tf})=1+$ (terms of higher degree), we get in particular 
for the
relative tangent bundle (assumed w.l.o.g. to 
 be a Hermitian bundle)
\begin{equation}\label{TdDev}
\ar\Td_{g_t}(\mtr{Tf})=\ar c^{\rm top}_{t}(\mtr N)^{-1}+\sum_{j=0}^{d-\rk N+1}
t^{-\rk N-j}\hat p_j(t)
\end{equation}
with classes $\hat p_j$ of degree $j$ which have a Taylor expansion in $t$ and 
$\hat
p_j(0)=0$.  Let $\deg_Y\alpha$ denote the degree of a Chow class $\alpha$ and 
define
$\deg_t t^k:=k$ for $k\in{\bf Z}$. Then any component $\alpha_t$ of the power 
series
$\ar\Td_{g_t}(\mtr{Tf})$ satisfies
$$
(\deg_Y+\deg_t)\alpha_t\geq -\rk N
$$
and equality is achieved precisely for the summand $\ar c^{\rm
top}_{t}(\mtr N)^{-1}$. Furthermore, by Lemma \ref{chgamma}
\begin{equation}\label{chDev}
\ar \ch_{g_t}(x)=\prod_{j=1}^k \left(\ar c_{q_j,t}(\mtr E^j)
+\ar s_j(t)\right)
\end{equation}
with classes $\ar s_j(t)$ of degree larger than $q_j$ which have a Taylor 
expansion in
$t$ and
$\ar s_j(0)=0$. Hence, any component $\alpha_t$ of the power series $\ar
\ch_{g_t}(x)$ satisfies
$$
(\deg_Y+\deg_t)\alpha_t\geq d+1
$$
and equality holds iff $\alpha_t$ is in the $\prod_{j=1}^k \ar c_{q_j,t}(\mtr
E^j)$-part. Hence any component $\alpha_t$ of $\ar\Td_{g_t}(\mtr{Tf})\ar
\ch_{g_t}(x)$ satisfies
\begin{equation}\label{deg1}
(\deg_Y+\deg_t)\alpha_t\geq d-\rk N+1\,\,.
\end{equation}
In particular the product has no singular terms in $t$, as $\deg_Y\leq d-\rk 
N+1$. In
other words, by multiplying formulae (\ref{TdDev}) and (\ref{chDev}) one obtains
$$
\ar\Td_{g_t}(\mtr{Tf})\ar \ch_{g_t}(x) =\ar c^{\rm top}_{t}(\mtr 
N)^{-1}\prod_{j=1}^k
\ar c_{q_j,t}(\mtr E^j) +O(t)\,\,,
$$
and the first summand on the right hand side has $\deg_Y=d-\rk N+1$, thus it is
constant in
$t$. We get formula (\ref{atregul}) by setting $t=1$. 
Now choose $\epsilon>0$ such that the fixed point set of $g_t$ does not vary on
$t\in]0,\epsilon[$.
To prove the second formula, we proceed similarly and we formally split 
$TX$ as a topological vector bundle into line bundles with first Chern class 
$x_\nu$,
acted upon by
$K$ with an angle
$\theta_\nu$. 
 The formulae for the
Lerch zeta function in \cite[p. 108]{Koehler2} or in \cite[Th. 7.10]{B1} show 
that the
$R$-class is given by
$$
R_{g_t}(TM)=-\sum_{\theta\neq0}\frac1{x_\nu+i
t\theta_\nu}\left(-2\Gamma'(1)+2\log|t\theta_\nu|+\log(1+\frac{x_\nu}{i
t\theta_\nu})\right)+O(1)
$$
for $t\to0$. Note that the singular term is of the form
$\alpha_1(t)\log |t|+\alpha_2(t)$ with
$$
(\deg_Y+\deg_t)\alpha_\mu(t)\geq -1
$$
($\mu=1,2$). As $(c^{\rm top}_{tK}(N)^{-1})'=c^{\rm top}_{tK}(N)^{-1}
\sum_{\theta\neq0}\frac{-1}{x_\nu+i t\theta_\nu}$ by definition, one obtains
\begin{eqnarray*}
\lefteqn{
\Td_{g_t}(TM)R_{g_t}(TM)\ch_{g_t}(x)
}\\&=&
\left((c^{\rm top}_{tK}(N)^{-1})'\log (t^2)+\frac{r_{tK}(N)}{c^{\rm
top}_{{tK}}(N)}\right)
\prod_j c_{q_j,tK}(E^j)+O(t)
\end{eqnarray*}
where the first term on the right hand side is again independent of $t$ (except 
the $\log (t^2)$).
\endProof
Note that the arithmetic Euler characteristic 
has a Taylor expansion in $t$.
Thus we get using Proposition \ref{regul} and Lemma \ref{rational}
\begin{cor} The equivariant analytic torsion of $x$ on the $d$-dimensional
K\"ahler manifold
$M$ has an asymptotic expansion
for
$t\to 0$
$$
T_{g_t}(M,\prod_{j} \gamma^{q_j}(\mtr E^j-\rk E^j))
=\log (t^2)\int_{M_K}\prod_j c_{q_j,K}(\mtr E^j) (c^{\rm top}_K(\mtr
N)^{-1})'+C_0+O(t\log |t|)
$$
with $C_0\in{\bf C}$.
\end{cor}
A more general version of this corollary is a consequence of \cite{BG} (see the 
next
section). 
We now combine our results with the (non-equivariant) arithmetic Riemann-Roch 
theorem.
We compute
\begin{eqnarray*}\lefteqn{
\ar{\deg}(f_{*}(\prod_{j}\ar{c}_{q_{j}}(\mtr{E}^j)))
}\\&=&
\ar{\deg}(f_{*}(\ar{\Td}(\mtr{Tf})\prod_{j}\ar{c}_{q_{j}}(\mtr{E}^j)))
=\ar{\deg}(f_*(\ar{\Td}(\mtr{Tf})\ar{\ch}(\prod_{j}\gamma^{q_j}(\mtr{E}^{j}-
\rk(E^j)))))
\\&=&
\ar{\deg}(f_{*}(\prod_{j}\gamma^{q_j}(\mtr{E}^{j}-
\rk(E^j))))-
{1\over 2}T(\prod_{j}\gamma^{q_j}(\mtr{E}^{j}-
\rk(E^j)))
\\&=&
L_{\torus}(\prod_{j}\gamma^{q_j}(\mtr{E}^{j}-
\rk(E^j)))-{1\over 2}T(\prod_{j}\gamma^{q_j}(\mtr{E}^{j}-
\rk(E^j)))
\\&&+
A_{\torus}(\prod_{j}\gamma^{q_j}(\mtr{E}^{j}-
\rk(E^j)))
\\&=&
L_\torus(\prod_{j}\gamma^{q_j}(\mtr{E}^{j}-
\rk(E^j)))
-{1\over 2}T(\prod_{j}\gamma^{q_j}(\mtr{E}^{j}-
\rk(E^j)))
\\&&
+\ar{\deg}(f_{*}^\torus({\prod_j\ar{c}_{q_j,t}(\mtr{E}^j)
\over \ar{c}_{t}^{\rm top}(\mtr{N})}))
\end{eqnarray*}
The first equality is justified by the fact that the $0$-degree part of the Todd 
class 
is $1$; the second one is \ref{cruxlem}; the third one is 
justified by the arithmetic Riemann-Roch theorem (\cite[eq. (1)]{GS8});  the 
fourth 
one is justified by (\ref{BC1}) and the last one by the last proposition. 
Finally, we get the following residue 
formula:
\begin{eqnarray}\lefteqn{
\ar\deg f_{*}(\prod_{j}\ar{c}_{q_{j}}(\mtr{E}^j))=L_{\torus}
(\prod_{j}\gamma^{q_j}(\mtr{E}^{j}-
\rk(E^j)))
}\nonumber\\&&
-\frac1{2}T(Y({\Bbb C}),\prod_{j}\gamma^{q_j}(\mtr{E}^{j}-
\rk(E^j)))+\ar\deg
f_*^{\torus}({\prod_{j}\ar{c}_{q_j,t}(\mtr{E}^j)\over 
\ar{c}^{\rm top}_{t}(\mtr{N})}).
\label{ResFV}
\end{eqnarray}
In particular, if $\mtr{L}$ is a $\torus$-equivariant line 
bundle on $Y$, one obtains the following formula for the height $h_{Y}(\mtr{L})$ 
of $Y$ relatively to $\mtr{L}$:
\begin{eqnarray*}
h_{\mtr{L}}(Y):=\ar{\deg}(\ar{c}_{1}(\mtr{L})^{d+1})&=&L_{\torus}((\mtr{L}-1)^{d
+1})
\\&&-
\frac1{2}T(Y({\Bbb C}),(\mtr{L}-1)^{d+1})+
\ar\deg f_*^{\torus}({\ar{c}_{1,t}(\mtr{L})^{d+1}\over 
\ar{c}^{\rm top}_{t}(\mtr{N})}).
\end{eqnarray*}
In our final residue formula, we shall use results of 
Bismut-Goette to give a formula for the term $L_{\torus}(\cdot)-
\frac1{2}T(Y({\Bbb C}),\cdot)$. 
Notice however that the last identity already implies 
the following corollary:
\begin{cor}
Let $Y$ be an arithmetic variety endowed with a $\torus$-action. Suppose 
that $\mtr{L}, \mtr{L}\pr$ are $\torus$-equivariant hermitian line bundles.
 If there is an equivariant isometry $\mtr{L}_{Y_\torus}\simeq \mtr{L}\pr_{Y_\torus}$ over 
 $Y_\torus$ and an equivariant (holomorphic) isometry 
 $\mtr{L}_{\bf C}\simeq\mtr{L}\pr_{\bf C}$ over $Y_{\bf C}$ then 
 $h_{\mtr{L}}(Y)=h_{\mtr{L}\pr}(Y)$. 
 \label{impcor}
\end{cor}

\subsection{The limit of the equivariant torsion}

Let $K'$ denote any nonzero multiple of $K$. The vector field
$K'$ is hamiltonian with respect to the K\"ahler form as the action on
$M$ factors through a projective space. Let $M_{K'}=M_K$ denote the fixed point 
set with
respect to the action of $K'$.
For any equivariant holomorphic Hermitian vector bundle $\mtr F$ we denote by
$\mu^F(K')\in\Gamma(M,\End(F))$ the section given by the action of the 
difference of
the Lie derivative and the covariant derivative
$L^F_{K'}-\nabla^F_{K'}$ on
$F$. Set as in \cite[ch. 7]{BGV}
$$
\Td_{K'}(\mtr {TM}):=\Td(-\frac{\Omega^{TM}}{2\pi i}+\mu^{TM}(K'))
\in{\frak A}(M)
$$
and
$$
\ch_{K'}(\mtr F):=\Tr \exp(-\frac{\Omega^F}{2\pi i}+\mu^F(K'))\in{\frak 
A}(M)\,\,.
$$
The Chern class $c_{q,K'}(\mtr F)$ for $0\leq q\leq \rk F$ is defined as the 
part of
total degree
$\deg_Y+\deg_t=q$ of
$$
\det(\frac{-\Omega^F}{2\pi i}+t\mu^F(K')+{\rm Id})
$$
at $t=1$, thus $c_{q,K'}(\mtr F)=c_q(-\Omega^F/2\pi i+\mu^F(K'))$.
Let $K^{\prime *}\in T^*_{\bf R}M$ denote the 1-form dual to $K'$ via
the metric on
$T_{\bf R}M$, hence $\iota_{K'} K^{\prime *}=\|K'\|^2$ is the norm square in 
$T_{\bf R}M$.
Set
$d_{K'}K^{\prime *}:=(d-2\pi i \iota_{K'})K^{\prime *}$ and define
$$
s_{K'}(u):=\frac{-\omega^{TM}}{2\pi u}\exp(\frac{d_{K'}K^{\prime *}}{4\pi
iu})=\sum_{\nu=0}^{d-1}
\frac{-\omega^{TM} (dK^{\prime *})^\nu}{2\pi u (4\pi i u)^\nu \nu!} 
e^{-\|K'\|^2/2 u}\,\,.
$$
For a smooth differential form $\eta$ it is shown in \cite{BG} (see also 
\cite[section
C,D]{Bmathann}) that the following integrals are well-defined:
$$
A_{K'}(\eta)(s):=\frac1{\Gamma(s)}\int_0^1\int_M\eta s_{K'}(u)
u^{s-1}du
$$
for Re $s>1$ and
$$
B_{K'}(\eta)(s):=\frac1{\Gamma(s)}\int_1^\infty\int_M\eta s_{K'}(u)
u^{s-1}du
$$
for Re $s<1$. Also it is shown in \cite{BG} (compare \cite[Proof of
theorem 7]{Bmathann}) that $s\mapsto A_{K'}(\eta)(s)$ has a meromorphic 
extension to $\bf C$
which is holomorphic at $s=0$ and that 
\begin{eqnarray*}\lefteqn{
A_{ K'}(\eta)'(0)+B_{ K'}(\eta)'(0)=
\int_1^\infty \int_M \eta 
s_{ K'}(u) \frac{du}{u}
}
\\&&+
\int_0^1 \int_M \eta \left(
s_{ K'}(u)
+\left(
\frac{\omega^{TM}}{2\pi u} c^{\rm top}_{K'}(\mtr N)^{-1}-(c^{\rm top}_{K'}(\mtr 
N)^{-1})'
\right)\delta_{M_{K'}}
\right)\frac{du}{u}
\\&&
+\int_{M_{K'}}\eta
\left(
\frac{\omega^{TM}}{2\pi} c^{\rm top}_{K'}(\mtr N)^{-1}-\Gamma'(1)(c^{\rm
top}_{K'}(\mtr N)^{-1})'
\right)
\end{eqnarray*}
for the derivatives $A_{ K'}(\eta)'$, $B_{ K'}(\eta)'$ of $A_{ K'}(\eta)$, $B_{
K'}(\eta)$ with respect to $s$; also
$$
A_{ K'}(\eta)(0)+B_{ K'}(\eta)(0)=\int_{M_{K'}}\eta (c^{\rm top}_{K'}(\mtr
N)^{-1})'\,\,.
$$
Define the current $S_{K'}(M,\omega^{TM})$ by the relation
$$
\int_M \eta S_{K'}(M,\omega^{TM}):=A_{K'}(\eta)'(0)+B_{K'}(\eta)'(0)\,\,.
$$
In particular, one notices
\begin{eqnarray*}
\lefteqn{
\int_M \eta S_{K'}(M,\omega^{TM})
=\lim_{a\to0^+}\Bigg[\int_a^\infty \int_M \eta 
s_{ K'}(u) \frac{du}{u}
}\\&&+
\int_a^1 \int_M \eta \left(
\frac{\omega^{TM}}{2\pi u} c^{\rm top}_{K'}(\mtr N)^{-1}-(c^{\rm top}_{K'}(\mtr 
N)^{-1})'
\right)\delta_{M_{K'}}
\frac{du}{u}
\\&&
+\int_{M_{K'}}\eta
\left(
\frac{\omega^{TM}}{2\pi} c^{\rm top}_{K'}(\mtr N)^{-1}-\Gamma'(1)(c^{\rm
top}_{K'}(\mtr N)^{-1})'
\right)\Bigg]
\\&=&\lim_{a\to0^+}\Bigg[
\int_M\eta \cdot 2i\omega^{TM}\frac{1-\exp(\frac{d_{K'}K^{\prime *}}{4\pi
ia})}{d_{K'}K^{\prime *}}
\\&&
+\int_{M_{K'}}\eta
\left(
\frac{\omega^{TM}}{2\pi a} c^{\rm top}_{K'}(\mtr N)^{-1}-(\Gamma'(1)+\log 
a)(c^{\rm
top}_{K'}(\mtr N)^{-1})'
\right)\Bigg]\,\,.
\end{eqnarray*}
By Lemma \ref{rational} and Proposition \ref{regul}, we already know that
$$
\lim_{t\to 0} \left(T_{g_t}(M,x)
-\log (t^2)\int_{M_K}\prod_j c_{q_j,K}(\mtr E^j) \cdot (c^{\rm top}_K(\mtr
N)^{-1})'\right)
$$
exists and equals $2 L_{\torus}(x)$. Now we shall compute this limit.
\begin{theor}\label{BismutGoette}
The limit of the equivariant analytic torsion of $x=\prod_{j} \gamma^{q_j}(\mtr
E^j-\rk E^j)$ associated to the action of $g_t$ for
$t\to 0$ is given by
\begin{eqnarray*}\lefteqn{
\lim_{t\to 0} \left(T_{g_t}(M,x)
-\log (t^2)\int_{M_K}\prod_j c_{q_j,K}(\mtr E^j) \cdot (c^{\rm top}_K(\mtr
N)^{-1})'\right) }\\&
=& T(M, x)+\int_M \prod_j c_{q_j,K}(\mtr E^j) \cdot S_K(M,\omega^{TM})
\end{eqnarray*}
\end{theor}
\beginProof
Let $I_{K'}$ denote the additive equivariant characteristic class which is given 
for a
line bundle $L$ as follows: If $K'$ acts at the fixed point $p$ by an angle
$\theta\in{\bf R}$ on
$L$, than
$$
I_{K'}(L)_{|p}:=\sum_{k\neq 0}\frac{\log(1+\frac\theta{2\pi 
k})}{c_1(L)+i\theta+2 k\pi
i}\,\,.
$$
The main result of \cite{BG} (announced in \cite{BG1}) implies that for 
$t\in{\bf
R}\setminus\{0\}$,
$t$ sufficiently small, there is a power series $T_t$ in $t$ with $T_0=T(M,x)$ 
such
that
\begin{eqnarray*}
T_{g_t}(M,x)-T_t&=&\int_M \Td_{tK}(\mtr{TM})\ch_{tK}(x)
S_{t K}(M,\omega^{TM})
\\&&-\int_{M_{g}}\Td_{g_t}(TM)\ch_{g_t}(x) I_{t K}(N_{M_{g}/M})\,\,.
\end{eqnarray*}
For $t\to 0$, both $I_{t K}(N_{M_{g}/M})\to 0$ and 
$\Td_{g_t}(\mtr{TM})\ch_{g_t}(x)\to
0$ (by eq. (\ref{deg1})), thus the last summand vanishes.

As in equation (\ref{chDev}) $\ch_{t K}(x)=\prod_j c_{q_j,tK}(\mtr 
E^j)+\tilde\eta(t)$
with a form $\tilde\eta$ such that $(\deg_Y+\deg_t)\tilde\eta(t)>d+1$. Thus
$\Td_{tK}(\mtr{TM})\ch_{t K}(x)=\prod_j c_{q_j,tK}(\mtr E^j)+\eta(t)$ with 
$(\deg_Y+\deg_t)\eta(t)>d+1$. Also 
$(\deg_Y+\deg_t)s_{t K}(t^2 u)=-1$, hence we observe that
$$
\int_M\Td_{t K}(\mtr{TM}) \ch_{t K}(x) s_{tK}(t^2u)
=\int_M \left(\prod_j c_{q_j,K}(\mtr E^j) s_K(u)+\eta(t) 
s_{tK}(t^2u)\right)\,\,.
$$
Let $\tilde\eta(t)$ denote the form obtained from $\eta(t)$ by multiplying the 
degree
$\deg_Y=j$ part with $t^{-j-1}$ for $0\leq j\leq d$. By making the change of 
variable
from
$u$ to
$t^2 u$ we get
$$
(A_{tK}+B_{tK})(\Td_{tK}(\mtr{TM})\ch_{tK}(x))(s)
=
t^{2 s}(A_{K}+B_{K})(\prod_j c_{q_j,K}(\mtr E^j)+\tilde\eta(t))(s)\,\,.
$$
Thus we find
\begin{eqnarray*}\lefteqn{
\int_M \Td_{tK}(\mtr{TM})\ch_{tK}(x)
S_{t K}(M,\omega^{TM})
}\\
&=&\log (t^2) \cdot(A_{K}+B_{K})(\prod_j c_{q_j,K}(\mtr E^j))(0)
\\&&
+(A_{K}+B_{K})(\prod_j c_{q_j,K}(\mtr E^j))'(0)+O(t\log (t^2))
\end{eqnarray*}
which implies  the statement of the theorem.
\endProof

\subsection{The residue formula}

By combining equation (\ref{ResFV}) and Theorem \ref{BismutGoette}, we obtain 
the following formula. 
Recall that $\torus$ is the one-dimensional diagonalisable torus over 
$\Spec{\Bbb
Z}$,  that $f:Y\ra\Spec\,{\Bbb Z}$ is a flat, $\torus$-projective morphism and 
that the fixed scheme $f^\torus:Y_{\torus}\ra\Spec\,{\Bbb Z}$ is assumed to be 
flat over $\Spec\,{\Bbb Z}$. We let 
$d+1$ be the absolute dimension of $Y$. 
We choose $\torus$-equivariant Hermitian bundles $\mtr{E}^{j}$ on $Y$ and 
positive integers $q_{j}$ such that $\sum_{j}q_{j}=d+1$. We get:
\begin{theor}
\begin{eqnarray*}
\lefteqn{
\ar{\rm deg}\bigg(f_{*}\big(\prod_{j}\ar{c}_{q_{j}}(\mtr{E}^j)\big)\bigg)=
\ar{\deg}\bigg(f_*^{\torus}\big({\prod_{j}\ar{c}_{q_j,t}(\mtr{E}^j)\over 
\ar{c}^{\rm top}_{t}(\mtr{N})}\big)\bigg)
}\\&&
+{1\over 2}\int_{Y({\Bbb C})}\prod_{j}{c}_{q_j,K}(\mtr{E}^j)\cdot S_{K}(Y({\Bbb 
C}),
\omega^{TY({\Bbb C})})
-{1\over 2}\int_{Y_{\torus}({\Bbb C})}\prod_{j}{c}_{q_j,K}(E^j)
\cdot {r_{K}(N) \over c_{K}^{\rm top}(N)}
\,\,.
\end{eqnarray*}
\label{residue}
\end{theor}

\par{\bf Example.} Assume that the fixed point scheme is flat of Krull dimension 
1.
The normal bundle to $Y_\torus$ splits as
$N=\bigoplus_{n\in{\bf Z}} N_n$. Thus 
$$
\hat c^{{\rm top}}_{t}(\mtr N)^{-1}=\frac1{\prod_n (2\pi it n)^{\rk N_n}}\left(
1-\sum_{n\in{\bf Z}} \frac{\hat c_1(\mtr{N}_n)}{2\pi it n}\right)
$$
by equation (\ref{cInv}). Also, at
a given point $p\in M_K$ the tangent space decomposes as
$TM_{|p}=\bigoplus TM_{\theta_\nu}$, where $K$ acts with angle
$\theta_\nu$ on $TM_{\theta_\nu}$.
Then
$$
\frac{r_K(N)}{c^{\rm top}_K(N)}_{|p}=\frac1{\prod_\theta i\theta}\sum_\theta
\frac{2\Gamma'(1)-2\log|\theta|}{i\theta}
$$
where the $\theta$ are counted with their multiplicity. Furthermore, in this 
case
\begin{eqnarray*}\lefteqn{
\int_M \eta S_{K}(M,\omega^{TM})=\lim_{a\to0^+}\Bigg[
\int_M \eta \cdot 2i\omega^{TM}\frac{1-\exp(\frac{d_{K}K^{
*}}{4\pi ia})}{d_{K}K^{*}}
}\\&&-(\Gamma'(1)+\log a)\int_{M_{K}}\eta\cdot
(c^{\rm
top}_{K}(\mtr N)^{-1})'
\Bigg]
\\&=&\int_M\eta\cdot 2i\omega^{TM}\left(
\frac1{d_K K^*}+\frac{(dK^*)^{d-1}}{(2\pi i\|K\|)^d}
\right)
\\&&-\lim_{a\to0^+}\Bigg[
\int_M\eta\cdot 
\frac{2i\omega^{TM}(dK^*)^{d-1}}{(2\pi i\|K\|)^d}(1-e^{-\frac1{2a}\|K\|^2})
\\&&+(\Gamma'(1)+\log a)\int_{M_{K}}\eta\cdot
(c^{\rm
top}_{K}(\mtr N)^{-1})'\Bigg]\,\,.
\end{eqnarray*}
Now consider a line bundle $\mtr {\cal L}$, splitting as 
$\bigoplus_k \mtr{\cal L}_k$ on
the fixed point scheme (where the ${\cal L}_{k}$ are locally free of 
rank $\leq 1$). We find 
\begin{eqnarray*}
\ar c_{1,t}(\mtr {\cal L})^{d+1}&=&\sum_k(\ar c_1(\mtr
{\cal L}_k)+2\pi i t k \,\rk {\cal L}_k)^{d+1}
\\&=&
\sum_k \left((2\pi i t k)^{d+1}\rk
{\cal L}_k+(d+1)(2\pi i t k)^d
\hat c_1(\mtr {\cal L}_k)\right)\,\,,
\end{eqnarray*}
thus
\begin{eqnarray*}
\lefteqn{
\widehat{\rm deg}\,f_*^{\torus}\left(\frac{\ar c_{1,t}(\mtr {\cal L})^{d+1}} 
{\ar
c^{\rm top}_t(\mtr N)} 
\right)
}\\&=&\ar\deg f_*^{\torus}\sum_k
\frac{k^d}{\prod_n n^{\rk N_n}}\left(
(d+1)\ar c_1(\mtr {\cal L}_k) -k\, \rk {\cal L}_k \cdot\sum_{n\in{\bf Z}} 
\frac{\ar
c_1(\mtr{N}_n)}{n}\right).
\end{eqnarray*}
Now notice that at a given fixed point $p$ over
$\bf C$ all but one ${\cal L}_{k,\bf C}$ vanish and set
$\phi_p:=2 \pi k$. We compute
$$
-\frac1{2}\sum_{p\in M_K} \frac{ c_{1,K}(L)^{d+1}} {
c^{\rm top}_K(N)} r_K(N)
=
\sum_{p\in M_K}
\frac{\phi_p^{d+1}}{\prod_\theta\theta}\sum_\theta\frac{-\Gamma'(1)+
\log|\theta|}{\theta}
$$
and
$$
-\frac1{2}(\Gamma'(1)+\log a)\sum_{p\in M_K}
c_{1,K}(L)^{d+1} (c^{\rm
top}_{K}(\mtr N)^{-1})'
=\sum_{p\in M_K}
\frac{\phi_p^{d+1}}{\prod_\theta\theta}\sum_\theta\frac{\Gamma'(1)+
\log a}{2\theta}\,\,.
$$
Hence we finally get
\begin{eqnarray*}
\ar\deg f_*\ar c_1(\mtr {\cal L})^{d+1}
&=&\ar\deg f_*^{\torus}\sum_k
\frac{k^d}{\prod_n n^{\rk N_n}}\left(
(d+1)\ar c_1(\mtr {\cal L}_k) -k\, \rk {\cal L}_k \cdot\sum_{n\in{\bf Z}} 
\frac{\ar
c_1(\mtr{N}_n)}{n}\right)
\\&&
+\sum_{p\in M_K}
\frac{\phi_p^{d+1}}{\prod_\theta\theta}\sum_\theta\frac{-\Gamma'(1)+
\log(\theta^2)}{2\theta}
\\&&+
\int_M c_{1,K}(L)^{d+1}\cdot i\omega^{TM}\left(
\frac1{d_K K^*}+\frac{(dK^*)^{d-1}}{(2\pi i\|K\|)^d}
\right)
\\&&-\lim_{a\to0^+}\Bigg[
\int_M \mu^L(K)^{d+1}\cdot 
\frac{i\omega^{TM}(dK^*)^{d-1}}{(2\pi i\|K\|)^d}(1-e^{-\frac1{2a}\|K\|^2})
\\&&-\log a
\sum_{p\in M_K}
\frac{\phi_p^{d+1}}{\prod_\theta\theta}\sum_\theta\frac{1}{2\theta}\Bigg]\,\,.
\end{eqnarray*}

\section{Appendix: a conjectural relative fixed point 
formula in Arakelov theory}

Since the first part of this series of articles was written, Xiaonan Ma 
defined in \cite{Ma2} higher analogs of the equivariant analytic torsion and 
proved curvature and anomaly formulae for it 
(in the case of fibrations by 
tori, this had already been done in \cite{Koehler4}). 
Once such torsion forms are available, one can 
formulate a conjectural fixed point formula, which fully 
generalizes \cite[Th. 4.4]{KR2} to the relative setting. 
Let $G$ be a 
compact Lie group and let $M$ and 
$M\pr$ be complex manifolds on which $G$ acts by holomorphic 
automorphisms. Let 
$f:M\ra M\pr$ be a smooth $G$-equivariant morphism of complex manifolds. 
Let 
${\omega}^{TM}$ be a $G$-invariant K\"ahler  metric on $M$ (a K\"ahler fibration
structure would in fact be sufficient). Let $\mtr{E}$ be a 
$G$-equivariant 
Hermitian holomorphic  vector bundle on $M$ and suppose for simplicity that 
$R^{k}f_{*}E=0$ for $k>0$. Now 
let $g$ be the automorphism corresponding to 
some element of $G$. The {\bf equivariant higher analytic torsion}  
$T_{g}(f,\mtr{E})$ is a certain element of $\widetilde{\frak A}(M_{g}\pr)$, 
which 
satisfies the curvature formula
$$
\dbd T_{g}(f,\mtr{E})=\ch_{g}(f_{*}\mtr{E})-
\int_{M_{g}/M_{g}\pr}\Td_{g}(\mtr{Tf})\ch_{g}(\mtr{E}).
$$
where $\mtr{Tf}$ is endowed with the metric induced by ${\omega}^{TM}$. 
The term in degree zero of $T_{g}(f,\mtr{E})$ is 
the equivariant analytic torsion $T_{g}(f^{-1}(x),\mtr{E}|_{f^{-1}(x)})$  
of the restriction of $\mtr{E}$ to the 
fiber of $f$ over $x\in M\pr_{g}$.\\
Now let $Y$, $B$ be $\mn$-equivariant arithmetic varieties over 
some fixed arithmetic ring $D$ and let $f:Y\ra B$ be 
a map over $D$, which is flat, $\mn$-projective and 
smooth over the complex 
numbers. Fix an $\mn({\bf C})$-invariant 
K\"ahler metric on $Y({\Bbb C})$. 
If $\mtr{E}$ is an $f$-acyclic (meaning that $R^k f_{*}E=0$ if $k>0$) 
$\mn$-equivariant Hermitian bundle 
on $Y$, let $f_{*}\mtr{E}$ be the direct image sheaf 
(which is locally free), 
endowed with its 
natural equivariant structure and $L_{2}$-metric.  
Consider the rule which associates the element 
$f_{*}\mtr{E}-T_{g}(f,\mtr{E})$ of 
$\ar K^{\mn}_{0}({B})$ to every $f$-acyclic equivariant 
Hermitian bundle $\mtr{E}$ 
and the element 
$\int_{Y({\Bbb C})_{g}/B({\Bbb C})_g}\Td_{g}(\mtr{Tf})\eta\in
\widetilde{\frak A}(B\lmn)$ to every 
$\eta\in\widetilde{\frak A}(Y_{\mn})$. The proof 
of the following proposition is then similar to 
the proof of \cite[Prop. 4.3]{KR2}. 
\begin{prop}
The above rule induces a group homomorphism 
$f_{*}:\ar K^{\mn}_{0}(Y)\ra \ar K^{\mn}_{0}(B)$.
\end{prop} 
Let $\cal R$ be the ring appearing 
in the statement of \cite[Th. 4.4]{KR2}. We are ready 
to formulate the following conjecture. 
\begin{conj}
Let
$$
M(f)=(\lambda_{-1}(\mtr{N}_{Y/Y_{\mu_{n}}}^{\vee}))^{-1}
\lambda_{-1}(f^{*}\mtr{N}_{B/B_{\mu_{n}}}^{\vee})
(1-R_{g}(N_{Y/Y_{\mu_{n}}})+R_{g}(f^{*}N_{B/B_\mn})).
$$
The diagram
$$
\matrix{
\ar{K}^{\mu_{n}}_{0}(Y)&
\stackrel{M(f).\rho(\cdot)}{\longrightarrow} 
&\ar{K}^{\mu_{n}}_{0}(Y_{\mu_{n}})\otimes_{R(\mu_{n})}{\cal R}\cr
   \downarrow\ f_{*}    &        &\downarrow\ f^{\mu_{n}}_{*}\cr
\ar{K}^{\mu_{n}}_{0}({B})&
 \stackrel{\rho(\cdot)}{\longrightarrow} 
&\ar{K}^{\mu_{n}}_{0}({B_\mn})\otimes_{R(\mu_{n})}{\cal R}\cr
}
$$
commutes. 
\label{PMainThp}
\end{conj}
About this conjecture, we make the following comments:

(a) One can carry over the principle of the proof of \cite[Th. 4.4]{KR2}
to prove this conjecture, provided a generalization of 
the immersion formula \cite[Th. 0.1]{Bimm} is available (which 
is not the case at the moment). We shall however not 
go into the details of this argument.

(b) Without formal proof again, we notice that 
the conjecture holds, if $B_\mn({\Bbb C})$ has dimension 
$0$. In that case the torsion forms are not 
necessary to define the direct image and the 
proof of \cite[Th. 4.4]{KR2} pulls through 
altogether.

\end{document}